\documentclass[11pt]{article}   

\usepackage{amssymb}
\usepackage{amsfonts}
\usepackage{latexsym}
\usepackage{amsthm}
\usepackage{enumerate}
\usepackage{epsfig}
\usepackage{graphicx}
\usepackage{color}
\usepackage{float}
\usepackage{subfigure}
\usepackage{amsmath}
\usepackage{makeidx}
\usepackage{fancyhdr}
\usepackage{lastpage}
\usepackage{url}
\newcommand{\Mod}[1]{\ (\text{mod}\ #1)}

\DeclareSymbolFont{AMSb}{U}{msb}{m}{n}  
\DeclareMathSymbol{\Sph}{\mathbin}{AMSb}{"53} \DeclareMathSymbol{\R}{\mathbin}{AMSb}{"52}
\DeclareMathSymbol{\T}{\mathbin}{AMSb}{"54} \DeclareMathSymbol{\Z}{\mathbin}{AMSb}{"5A}
\DeclareMathSymbol{\K}{\mathbin}{AMSb}{"4B}

\newtheorem{guess}{Theorem}[section]

\newtheorem{corollary}{Corollary}[section]
\newtheorem{lemma}{Lemma}[section]

\newtheorem{proposition}{Proposition}[section]

\def\qed{\hfill \rule{4pt}{7pt}}
\def\pf{\noindent {\it Proof. }}
\def\qed{\hfill \rule{4pt}{7pt}}

\newcounter{case}

\begin{document}  

\title{Constructing $2$-Arc-Transitive Covers of Hypercubes}

\author{Michael Giudici$^{1}$\footnote{Email: michael.giudici@research.uwa.edu.au.}, Cai Heng Li$^{2}$\footnote{Email: lich@sustc.edu.cn.}, Yian Xu$^{3}$\footnote{Email: yian.xu@research.uwa.edu.au.} \\\\
 $^{1}$ Centre for the Mathematics of Symmetry and Computation, \\
 School of Mathematics and Statistics,  \\
  The University of Western Australia \vspace{5pt} \\  
  $^{2}$ Department of Mathematics\\
   Southern University of Science and Technology\\
  $^{3}$ School of Mathematics\\
  Southeast University
}

 

\date{}
\maketitle

\begin{abstract}
We introduce the notion of a symmetric basis of a vector space equipped with a quadratic form, and provide a sufficient and necessary condition for the existence to such a basis. Symmetric bases are then used to study Cayley graphs of certain extraspecial $2$-groups of order $2^{2r+1}$ ($r\geq 1$), which are further shown to be normal Cayley graphs and $2$-arc-transitive covers of $2r$-dimensional hypercubes. 

\textit{Keywords} extraspecial 2-group \and symmetric basis \and quadratic form \and locally-primitive graph \and normal Cayley graph
\end{abstract}

\section{Introduction}
\label{intro}
Throughout this paper, all graphs are simple, connected and regular. Let $\Gamma$ be a graph with vertex set $V(\Gamma)$ and edge set $E(\Gamma)$. An {\it s-arc} of $\Gamma$ is a sequence $(v_{1}, \ldots, v_{s+1})$ of $s+1$ vertices such that for all $1\leq i\leq s$, $\{v_{i}, v_{i+1}\}$ is an edge in $\Gamma$ and $v_{i}\neq v_{i+2}$, and $\Gamma$ is said to be {\it s-arc-transitive} if the automorphism group of $\Gamma$ is transitive on the set of $s$-arcs. The study of $s$-arc-transitive graphs is motivated by a result of Tutte (1949), which says that there are no $s$-arc-transitive graphs of valency $3$ for $s\geq 6$. Later this result was extended by Weiss \cite{Weiss1981} saying that there are no $8$-arc-transitive graphs of valency at least 3. Thus analysing the $s$-arc-transitive graphs for $2\leq s\leq 7$ has become one of the central goals in algebraic graph theory, and the classification of some 2-arc-transitive graphs has been obtained. For example, the 2-arc-transitive circulants are classified in \cite{Alspach1996}; a complete classification of 2-arc-transitive dihedrants is given in \cite{Du2008}; and a class of 2-arc-transitive Cayley graphs of elementary abelian 2-groups is classified in \cite{ivanov1993finite}. 

A natural idea to investigate the 2-arc-transitive graphs is to study their quotient graphs. Let $\mathcal{P}$ be a partition of the vertex set $V(\Gamma)$. Define the {\it quotient graph} $\Gamma_{\mathcal{P}}$ of $\Gamma$ to be the graph with vertex set $\mathcal{P}$ and two parts $P, P'\in \mathcal{P}$ form an edge if and only if there is at least one edge in $\Gamma$ joining a vertex of $P$ and a vertex of $P'$. If $\mathcal{P}$ is {\it $G$-invariant} for some group $G$ of automorphisms of $\Gamma$, then the action of $G$ on $\Gamma$ induces an action of $G$ on $\Gamma_{\mathcal{P}}$. Let $N$ be a nontrivial normal subgroup of $G$ and $\mathcal{P}$ be the set of $N$-orbits in $V(\Gamma)$. The quotient graph $\Gamma_{\mathcal{P}}$ is said to be a {\it normal quotient} of $\Gamma$, denoted $\Gamma_{N}$. In general, the valency of $\Gamma_{N}$ divides the valency of $\Gamma$. If the valency of $\Gamma$ equals the valency of $\Gamma_{N}$, then $\Gamma$ is said to be a {\it cover} of $\Gamma_{N}$. It has been proved by Praeger \cite[Theorem 4.1]{praeger1993nan} that if $G$ is vertex-transitive and $2$-arc-transitive on $\Gamma$, and $N$ has more than two orbits in $V(\Gamma)$, then 
\begin{enumerate}
\item $G/N$ is $s$-arc transitive on $\Gamma_{N}$ and $G/N$ is faithful on $V(\Gamma_{N})$, 
\item $\Gamma$ is a cover of $\Gamma_{N}$, and 
\item $N$ is semiregular on $V(\Gamma)$. 
\end{enumerate}
We say that a permutation group is {\it quasiprimitive}  on a set $\Omega$ if every nontrivial normal subgroup of the permutation group is transitive on $\Omega$, and {\it primitive} if it acts transitively on $\Omega$ and preserves no nontrivial partition of $\Omega$.  A permutation group is said to be {\it bi-quasiprimitive} on $\Omega$ if 
\begin{itemize}
\item[](i) each nontrivial normal subgroup of the permutation group has at most two orbits on $\Omega$, and
\item[](ii) there exists a normal subgroup with two orbits on $\Omega$.
\end{itemize}
The structure of finite quasiprimitive permutation groups was investigated in \cite{praeger1993nan} and the types of quasiprimitive groups that are 2-arc-transitive on a graph were determined. Praeger studied bi-quasiprimitive groups in \cite{Praeger2003} and one specific class identified was previously studied in \cite{praeger1993bipartite}. One family of such bipartite bi-quasiprimitive graphs are the affine ones. A $2$-arc-transitive graph is said to be {\it affine}, if there is a vector space $N$ and a group $G$ of automorphisms of the graph such that $N\leq G\leq AGL(N)$ with $N$ regular on the vertices and $G$ acting transitively on the set of $2$-arcs. Table 1 in \cite{ivanov1993finite} classifies all affine bipartite $2$-arc-transitive graphs with the stabilizer of the bipartition of the vertices being primitive on each bipartition. 

Another interesting topic is to reconstruct $2$-arc-transitive covers of $2$-arc-transitive graphs. It is known that every finite regular graph has a $2$-arc-transitive cover \cite{babai1985arc}. In \cite{Du1998} Du, Malni\u{c} and Waller investigate the regular covers of complete graphs which are 2-arc-transitive, and they give a complete classification of all graphs whose group of covering transformations is either cyclic or isomorphic to $Z_{p}\times Z_{p}$ where $p$ is a prime and whose fibre-preserving subgroup of automorphisms acts 2-arc-transitively. In particular, two families of 2-arc-transitive graphs are obtained. After that, many more results related to the reconstruction of the 2-arc-transitive graphs have been obtained, see \cite{Du2005} for examples. 

The main subject of this paper is to construct a $2$-arc-transitive cover for one family of affine graphs, namely the hypercubes. Let $V=\mathbb{Z}_{2}^{d}$ be a $d$-dimensional vector space over the field $\mathbb{F}_{2}$, let $e_{1}, \ldots, e_{d}$ be a basis of $V$ and $\mathcal{B}=\{e_{1}, \ldots, e_{d}\}$. A {\it $d$-dimensional hypercube} is a Cayley graph $Q_{d}=Cay(V, \mathcal{B})$. It is known that $Q_{d}$ admits a regular group of automorphisms $Z_{2}^{d}$ and $Aut(Q_{d})=Z_{2}^{d}\rtimes S_{d}$, where $S_{d}=Aut_{\bf1}(Q_{d})$ and permutes $e_{1}, \ldots, e_{d}$ naturally (see \cite{spiga2009enumerating}). Thus the hypercubes are 2-arc-transitive affine graphs. Furthermore, it has been shown in \cite{ivanov1993finite} that $Q_{d}$ is bi-quasiprimitive if and only if $d=2$ or $d$ is odd. In this paper, we are interested in the even-dimensional hypercubes, in particular, we construct a $2$-arc-transitive cover for even-dimensional hypercubes. We also show that such a cover is a normal Cayley graph.

\medskip

Let $G$ be a finite group, and $S$ be a subset of $G$ such that $S$ does not contain the identity of $G$ and $S=S^{-1}=\{s^{-1}|s\in S\}$. We say that an element $g$ of $G$ is an {\it involution} if it has order $2$, that is, $g^{2}={\bf 1}$ and $g\neq \bf1$. The $\it{Cayley\ graph}$ $\Gamma=Cay(G, S)$ is defined to have vertex set $V(\Gamma)=G$, and edge set $E(\Gamma)=\{\{g, sg\}|s\in S\}$. It is well known that a graph is a Cayley graph if and only if its full automorphism group contains a subgroup acting regularly on the vertex set of the graph (see \cite{sabidussi1958class}). 

Let $\Gamma=Cay(G, S)$ be a Cayley graph for some group $G$ and $Aut(\Gamma)$ be the full automorphism group of $\Gamma$. For each $g\in G$, define a map $\hat{g}: G\to G$ by the right multiplication of $g$ on $G$ as below: 
$$
\hat{g}: x\to xg, \ for\ x\in G.
$$
Then $\hat{g}$ is an automorphism of $\Gamma$. It follows from the definition that the group $\hat{G}=\{\hat{g}\mid g\in G\}$ is a subgroup of $Aut(\Gamma)$ and acts regularly on $V(\Gamma)$. 
Following Xu \cite{xu1998automorphism}, we say that $\Gamma$ is a {\it normal\ Cayley\ graph\ for\ G} (or normal) if $\hat{G}\unlhd Aut(\Gamma)$, otherwise we say that $\Gamma$ is a {\it non}-{\it normal\ Cayley\ graph\ for\ G} (or non-normal). 

\medskip

Suppose that $V$ is a $d$-dimensional vector space with a nondegenerate quadratic form $Q$ where the associated bilinear form $B_{Q}$ is symmetric. Let $\mathcal{C}=\{v_{1}, v_{2}, \ldots, v_{d}\}$ be a basis of $V$. We say that $\mathcal{C}$ is {\it symmetric} if $Q(v_{i})=0$ and $B(v_{i}, v_{j})=1$ for all $i, j$ with $1\leq i<j\leq d$. In Section 3, we determine a necessary and sufficient condition for a vector space to have a symmetric basis.

\medskip

Let $G$ be an extraspecial $2$-group of order $2^{2r+1}$ with $r\geq 1$, that is, $|Z(G)|=2$ and $G/Z(G)\cong \mathbb{Z}_{2}^{2r}$. There are two extraspecial $2$-groups of each order, for which we will give more details in Section $3$. Let $\overline{S}=\{\overline{s_{1}}, \ldots, \overline{s_{2r}}\}$ be a symmetric basis of $G/Z(G)$, and for each $1\leq i\leq 2r$, let $s_{i}$ be a preimage of $\overline{s_{i}}$ in $G$. Notice that generally for a basis of $G/Z(G)$, the preimages of the basis elements are not necessary involutions in $G$. However in Section 3 we show that in the case where it is a symmetric basis the preimages of the basis elements are all involutions in $G$, which is crucial for the proofs of the main results.  Note that $\Sigma=Cay(G/Z(G), \overline{S})$ is a $2r$-dimensional hypercube. We will prove the following result. 

\medskip

\begin{guess}\label{ch4mainthem}
Let $G$ be an extraspecial $2$-group of order $2^{2r+1}$ with $r\geq 1$ such that $G/Z(G)$ has a symmetric basis $\{\overline{s_{1}}, \ldots, \overline{s_{2r}}\}$. Let $\Gamma=Cay(G, S)$ be a Cayley graph of $G$ with $S=\{s_{1}, \ldots, s_{2r}\}$. Then $\Gamma$ is a $2$-arc-transitive cover of some $2r$-dimensional hypercube $\Sigma$, and $\Gamma$ is a normal Cayley graph with $Aut(\Gamma)=G\rtimes S_{2r}$.
\end{guess}

\medskip

\section{Preliminaries}

Let $V$ be a $2r$-dimensional ($r\geq 1$) vector space over a field $\mathbb{F}_{q}$, where $q$ is a prime-power. Let $B$ be a bilinear form on $V$. We say that $B$ is {\it symmetric } if $B(u, v)=B(v, u)$ for all $u, v\in V$, and $B$ is {\it alternating} if $B(u, u)=0$ for all $u\in V$. The {\it radical} of $B$ is the subspace 
\begin{center}
$rad(B)=\{u\in V\mid B(u, v)=0$ for all $v\in V\}$,
\end{center}
and $B$ is said to be {\it nondegenerate} if rad$(B)=\{0\}$. Let $W$ be a subspace of $V$. Define 
\begin{equation*}
W^{\perp}=\{v\in V\mid B(w, v)=0 \text{ for all }w\in W\}
\end{equation*}
to be the {\it orthogonal complement} of W. It is known that if $B$ is nondegenerate, then $dim(W)+dim(W^{\perp})=dim(V)$.
A map $Q: V\to \mathbb{F}_{q}$ is a {\it quadratic form} on $V$ if the following two conditions are satisfied:
\begin{itemize}
\item[](i) $Q(\lambda u)=\lambda^{2}Q(u)$ for all $u\in V$ and $\lambda\in \mathbb{F}_{q}$, and
\item[](ii) the map $B_{Q}: V\times V\to \mathbb{F}_{q}$ defined by 
$$
B_{Q}(u, v)=Q(u+v)-Q(u)-Q(v)
$$
is a bilinear form.
\end{itemize} 
The bilinear form $B_{Q}$ is called the {\it associated bilinear form} of $Q$. A quadratic form $Q$ is said to be {\it nondegenerate} if and only if its associated bilinear form is nondegenerate. 

Let $u, v$ be two distinct vectors of $V$. We say that $\{u, v\}$ is a {\it hyperbolic pair} if $Q(u)=Q(v)=0$ and $B(u, v)=1$. 
    By \cite[Proposition 2.2.7]{burness2016classical}, when $Q$ is nondegenerate and $B=B_{Q}$ is symmetric , $V$ has the following two types of standard bases, in particular, the basis is {\it hyperbolic} in case (i), and {\it elliptic} in case (ii):
\begin{itemize}
\item[] (i). $\mathcal{B}=\{e_{1}, \ldots, e_{n/2}, f_{1}, \ldots, f_{n/2}\}$ where 
\begin{center}
$Q(e_{i})=Q(f_{i})=0$, $B(e_{i}, f_{j})=\delta_{ij}$ for all $i,j$;
\end{center} 
\item[] (ii). $\mathcal{B}=\{e_{1}, \ldots, e_{n/2-1}, f_{1}, \ldots, f_{n/2-1}, x, y\}$ where $Q(e_{i})=Q(f_{i})=0$,
\begin{center}
$B(e_{i}, x)=B(e_{i}, y)=B(f_{i}, x)=B(f_{i}, y)=0$, $B(e_{i}, f_{j})=\delta_{ij}$
\end{center} 
 for all $i,j$, $Q(x)=1$, $B(x,y)=1$ and $Q(y)=\zeta$ where $x^{2}+x+\zeta\in \mathbb{F}_{q}[x]$ is irreducible,
\end{itemize}
where 
\[ \delta_{ij}= \left\{
  \begin{array}{l l}
    1 & \quad \text{if $i=j$}\\
    0 & \quad \text{if $i\neq j$.}
  \end{array} \right.\]
If $V$ has a hyperbolic basis, then $Q$ is said to be a {\it hyperbolic} quadratic form (or hyperbolic in short). Similarly we say that $Q$ is {\it elliptic} when $V$ has an elliptic basis.  

Let $U$ be a subspace of $V$. We say that $U$ is {\it totally singular} if $Q(u)=0$ for all $u\in U$. The next result is a consequence of Witt's Lemma (see \cite{burness2016classical})
\begin{proposition}[{\cite[Page 38]{Burness2016}}]\label{ch4pro21}
Let $V$ be a $d$-dimensional vector space over the field $\mathbb{F}_{q}$ equipped with a nondegenerate quadratic form $Q$, and $U$ be a maximal totally singular subspace of $V$. Then 
$$
dimU=\frac{1}{2}dimV-\delta
$$
where $\delta=0$ if $Q$ is hyperbolic, and $\delta=1$ if $Q$ is elliptic.
\end{proposition}

\section{Symmetric Basis of a Vector Space}

Let $V$ be a $2r$-dimensional vector space over $\mathbb{F}_{2}$ with a nondegenerate quadratic form $Q: V\to \mathbb{F}_{2}$ and an associated symmetric  bilinear form $B=B_{Q}$. Suppose that $\mathcal{C}=\{v_{1}, \ldots, v_{2r}\}$ is a symmetric basis of $V$. Let $\mathcal{C}'=\{c_{1}, \ldots, c_{r-1}\}$ where for $1\leqslant  i\leqslant  r-1$, we have 
\begin{center}
$c_{i}=v_{2i-1}+v_{2i}+v_{2i+1}+v_{2i+2}$.
\end{center}
Let $U$ be the subspace of $V$ generated by $\mathcal{C}'$. Then $dim(U)=\frac{1}{2}dim(V)-1$. Also $U$ is totally singular as $Q(c_{i})=0$ and $B(c_{i}, c_{j})=0$ for all $1\leqslant  i, j\leqslant  r-1$. If $U$ is maximal subject to being totally singular in $V$, then $Q$ is elliptic. Otherwise $Q$ is hyperbolic, and there is a vector $v\in V\backslash U$ such that $W=\langle v, \mathcal{C}'\rangle$ is a maximal totally singular subspace of $V$. 

Let $\alpha\in V$ be a non-zero vector that is not in $\mathcal{C}'$. We may assume that $\alpha=\mu_{1}+\mu_{2}+\cdots+\mu_{t}$, where $\mu_{i}\in \mathcal{C}$ and $\mu_{i}\neq \mu_{j}$ for all $1\leqslant  i\neq j\leqslant  t$.

\medskip

\begin{lemma}\label{ch4lem0}
$Q(\alpha)=0$ if and only if $t\equiv 0$ or $1 \Mod 4$.
\end{lemma}

\smallskip

\pf  First suppose that $t=2k$ with $k\geqslant  0$ and let $u_{i}=\mu_{2i-1}+\mu_{2i}$, for $i=1,2,\ldots, k$. So 
\begin{align*}
& Q(\mu_{1}+\mu_{2}+\cdots+\mu_{2k})  \\
& = Q(u_{1}+u_{2}+\cdots+u_{k}) \\
& = Q(u_{1})+Q(u_{2}+\cdots+u_{k})+\sum\limits_{i=2}^{k}B(u_{1}, u_{i})  \\
& = Q(u_{1})+Q(u_{2})+Q(u_{3}+\cdots+u_{k})+\sum\limits_{i=2}^{k}B(u_{1}, u_{i})+\sum\limits_{i=3}^{k}B(u_{2}, u_{i})  \\
& = Q(u_{1})+Q(u_{2})+Q(u_{3})+Q(u_{4}+\cdots+u_{k})+\sum\limits_{i=2}^{k}B(u_{1}, u_{i})+\sum\limits_{i=3}^{k}B(u_{2}, u_{i})+\sum\limits_{i=4}^{k}B(u_{3}, u_{i})  \\
& = \sum\limits_{i=1}^{k}Q(u_{i})+\sum\limits_{i=1}^{k-1}\sum\limits_{j=i+1}^{k}B(u_{i}, u_{j}) .
\end{align*} 
For $1\leqslant  i<j\leqslant  k$ we have 
\begin{align*}
B(u_{i}, u_{j}) & =B(\mu_{2i-1}+\mu_{2i}, \mu_{2j-1}+\mu_{2j})\\
& =B(\mu_{2i-1}, \mu_{2j-1})+B(\mu_{2i-1}, \mu_{2j})+B(\mu_{2i}, \mu_{2j-1})+B(\mu_{2i}, \mu_{2j})\\
& =0, 
\end{align*}
and 
\begin{align*}
Q(u_{i}) =& Q(\mu_{2i-1}+\mu_{2i})\\
& =Q(\mu_{2i-1})+Q(\mu_{2i})+B(\mu_{2i-1}, \mu_{2i})\\
& =1.
\end{align*}
Thus $Q(\alpha)=0$ if and only if $k\equiv 0 \Mod 2$, that is, if and only if $t\equiv 0 \Mod 4$. 

Next suppose that $t$ is odd, and so $t=2k+1$ with $k\geqslant  0$. Then 
\begin{align*}
Q(\mu_{1}+\mu_{2}+\cdots+\mu_{2k}+\mu_{2k+1}) &=Q(\mu_{1}+\cdots+\mu_{2k})+Q(\mu_{2k+1})+\sum\limits_{i=1}^{2k}B(\mu_{i}, \mu_{2k+1})\\
& =Q(\mu_{1}+\cdots+\mu_{2k}).
\end{align*}
which implies that $Q(\alpha)=0$ if and only if  $t\equiv 1 \Mod 4$. \qed

\medskip

Before we introduce the sufficient and necessary conditions for $V$ to have a symmetric basis, we first give the following lemma. 

\begin{lemma}\label{ch4lem01}
Let $W$ be a nontrivial subspace of $V$ with dimension $d\equiv 2 \Mod 4$. Suppose that $W$ has a symmetric basis and there are three pairwise perpendicular hyperbolic pairs $\{a, b\}, \{c, d\}, \{g, h\}$ in $W^{\perp}\backslash W$. Let 
\begin{equation*}
W'=W\bot\langle a, b\rangle\bot\langle c, d\rangle\bot\langle g, h\rangle.
\end{equation*}
Then $W'$ has a symmetric basis.
\end{lemma}

\medskip

\pf 
Let $\mathcal{C}_{W}=\{w_{1}, \ldots, w_{d}\}$ be a symmetric basis of $W$. Let
\begin{center}
$u_{1}=a+c+d+\sum\limits_{i=1}^{d}w_{i}$;\\
$u_{2}=b+c+d+\sum\limits_{i=1}^{d}w_{i}$;\\
$u_{3}=c+g+h+\sum\limits_{i=1}^{d}w_{i}$;\\
$u_{4}=d+g+h+\sum\limits_{i=1}^{d}w_{i}$;\\
$u_{5}=g+a+b+\sum\limits_{i=1}^{d}w_{i}$;\\
$u_{6}=h+a+b+\sum\limits_{i=1}^{d}w_{i}$.
\end{center}
Since $d\equiv 2 \Mod 4$, Lemma~\ref{ch4lem0} implies $Q(\sum\limits_{i=1}^{d}w_{i})=1$. So $Q(u_{j})=0$ for all $j=1, 2,\ldots, 6$. For any $1\leqslant  j_{1}\neq j_{2}\leqslant  6$, one can check that
\begin{center}
$B(u_{j_{1}}, u_{j_{2}})=1$.
\end{center}
Also we have $B(u_{j}, w_{i})=1$ for all $1\leqslant  i\leqslant  d$ and $j=1, 2, \ldots, 6$. Therefore $\mathcal{C}_{W}\cup \{u_{1}, \ldots, u_{6}\}$ is a symmetric basis of $W'$. \qed 

\medskip

\begin{lemma}\label{ch4lem32}
Let $V$ be a vector space of dimension $2r$ with nondegenerate quadratic form $Q$ such that $r\geqslant  1$. 
\begin{itemize}
\item[](i). If $Q$ is hyperbolic and $r\equiv 0$ or $1 \Mod 4$, then $V$ has a symmetric basis;
\item[](ii). If $Q$ is elliptic and $r\equiv 2$ or $3\Mod 4$, then $V$ has a symmetric basis.
\end{itemize}
\end{lemma}

\pf We prove this by using induction on $r$. 

{\it (i)} Suppose that $Q$ is hyperbolic with $r\equiv 0$ or $1 \Mod 4$. When $r=1$, a hyperbolic basis of $V$ is also a symmetric basis of $V$. Now assume the lemma holds for all $r\leqslant  4\ell+1$ with $r\equiv 0$ or $1 \Mod 4$, where $\ell$ is a nonnegative integer. Note that we have seen that the lemma holds when $\ell=0$. Let $r'=4\ell$ and suppose that $r=r'+4=4(\ell+1)$. Let $W\leqslant  V$ be a subspace of dimension $2(r'+1)$ such that $Q$ is hyperbolic on $W$ and 
\begin{center}
$V=W\bot\langle a, b\rangle\bot\langle c, d\rangle\bot\langle g, h\rangle$.
\end{center} 
where $\{a, b\}, \{c, d\}, \{g, h\}$ are hyperbolic pairs in $V\backslash W$. By induction, $W$ has a symmetric basis $\mathcal{C}_{W}=\{w_{1}, w_{2}, \ldots, w_{2(r'+1)}\}$. Take $d=2(r'+1)$. Then $V$ has a symmetric basis $\mathcal{C}_{W}\cup \{u_{1}, \ldots, u_{6}\}$ as constructed in Lemma~\ref{ch4lem01}.  

Now suppose that $r=4(\ell+1)+1$. We may assume that 
\begin{center}
$V=W\bot\langle a, b\rangle\bot\langle c, d\rangle\bot\langle g, h\rangle\bot\langle x,y\rangle$,
\end{center}
where $\{a, b\}, \{c, d\}, \{g, h\}, \{x, y\}$ are hyperbolic. Clearly $V$ contains a subspace $U$ with a symmetric basis $\mathcal{C}_{W}\cup \{u_{1}, \ldots, u_{6}\}$ as above, and let 
\begin{center}
$\alpha=x+\sum\limits_{i=1}^{2r'+2}w_{i}+a+b+c+d+g+h$;\\
$\beta=y+\sum\limits_{i=1}^{2r'+2}w_{i}+a+b+c+d+g+h$.
\end{center}
Then $Q(\alpha)=Q(\beta)=0$ and $B(\alpha, \beta)=1$. Also for $1\leqslant  i\leqslant  2(r'+1)$ we have $B(\alpha, w_{i})=B(\beta, w_{i})=1$, and $B(\alpha, u_{j})=B(\beta, u_{j})=1$ for $1\leqslant  j\leqslant  6$. Hence $\mathcal{C}_{W}\cup \{u_{1}, \ldots, u_{6}\}\cup \{\alpha, \beta\}$ forms a symmetric basis of $V$.

\smallskip

{\it (ii)} Suppose that $Q$ is elliptic and $r\equiv 2$ or $3 \Mod 4$. When $r=2$, let $\{e_{1}, f_{1}, x, y\}$ be an elliptic basis of $V$. Let 
\begin{center}
$c_{1}=e_{1}$, $c_{2}=f_{1}$, $c_{3}=x+c_{1}+c_{2}$, $c_{4}=y+c_{1}+c_{2}$.
\end{center} 
Then $Q(c_{i})=0$ and $B(c_{i}, c_{j})=1$ for $1\leqslant  i<j\leqslant  4$. So $\{c_{1}, c_{2}, c_{3}, c_{4}\}$ is a symmetric basis of $V$. 

When $r=3$, let $\{e_{1}, e_{2}, f_{1}, f_{2}, x, y\}$ be an elliptic basis. Let $c_{1}, c_{2}, c_{3}, c_{4}$ be defined as above and let
\begin{center}
$c_{5}=e_{2}+c_{1}+c_{2}+c_{3}+c_{4}$, $c_{6}=f_{2}+c_{1}+c_{2}+c_{3}+c_{4}$.
\end{center} 
Then $\{c_{1}, c_{2}, c_{3}, c_{4}, c_{5}, c_{6}\}$ is a symmetric basis of $V$.

Now assume the lemma holds for all $r\leqslant  4\ell+3$ with $r\equiv 2$ or $3 \Mod 4$ where $\ell$ is a nonnegative integer. Note that we have seen that the lemma holds when $\ell=0$. Let $r'=4\ell+2$ and suppose that $r=r'+4=4(\ell+1)+2$. Let $W\leqslant  V$ be a subspace of dimension $2(r'+1)$ such that $Q$ is elliptic on $W$ and 
\begin{center}
$V=W\bot\langle a, b\rangle\bot\langle c, d\rangle\bot\langle g, h\rangle$,
\end{center} 
where $\{a, b\}, \{c, d\},\{g, h\}$ are hyperbolic. By our induction $W$ has a symmetric basis 
$$
\mathcal{C}_{W}=\{w_{1}, \ldots, w_{2r'+1}, w_{2r'+2}\}.
$$ 
Then by Lemma~\ref{ch4lem01}, when $r=4(\ell+1)+2$, $V$ contains a symmetric basis. 
When $r=4(\ell+1)+3$, let $x,y, \alpha, \beta$ be vectors of $V$ as defined in $(i)$. Then $\mathcal{C}_{W}\cup \{u_{1}, \ldots, u_{6}\}\cup \{\alpha, \beta\}$ forms a symmetric basis of $V$. \qed

Let $V$ be a $d$-dimensional vector space over field $\mathbb{F}_{q}$ equipped with a quadratic form $Q$ and $U$ be a $d$-dimensional vector spaces over field $\mathbb{F}_{q}$ equipped with a quadratic form $Q'$. An isometry from $V$ to $U$ is an invertible linear map $\sigma: V\to U$ such that
\begin{equation}\label{ch4equ21}
Q'(v^{\sigma})=Q(v)
\end{equation}
for all $v\in V$. Notice that (\ref{ch4equ21}) implies that  
\begin{equation*}
B_{Q'}(u^{\sigma}, v^{\sigma})=B_{Q}(u, v)
\end{equation*}
for all $u, v\in V$. If such an isometry exists, then both $U$ and $V$, and $Q$ and $Q'$ are said to be {\it isometric}. Let $W\subseteq V$ and $W'\subseteq U$. Then $W$ and $W'$ are isometric if the restrictions $Q\mid _{W}$ and $Q'\mid_{W'}$ are isometric. We say that $\sigma$ is an {\it isometry of $Q$} if $U=V$. The {\it isometry group} of $Q$ is the set of isometries of $Q$ under composition. Notice that the isometry group of $Q$ is a subgroup of the isometry group of the associated bilinear form $B_{Q}$. 

\begin{lemma}\label{ch4lem33}
Let $V$ be a vector space of dimension $2r$ with nondegenerate quadratic form $Q$ and a symmetric basis $\mathcal{C}=\{v_{1}, v_{2}, \ldots, v_{2r}\}$. 
\begin{itemize}
\item[](i). If $r\equiv 0$ or $1 \Mod 4$, then $Q$ is hyperbolic;
\item[](ii). If $r\equiv 2$ or $3 \Mod 4$, then $Q$ is elliptic. 
\end{itemize}
\end{lemma}
 
\pf Let $Q_{1}$ and $Q_{2}$ be quadratic forms on $V$ such that $V$ has a symmetric basis $\mathcal{C}=\{c_{1}, \ldots, c_{2r}\}$ with respect to $Q_{1}$ and a symmetric basis $\mathcal{C'}=\{c'_{1}, \ldots, c'_{2r}\}$ with respect to $Q_{2}$. Let $\sigma: V\to V$ be the linear map defined by
\begin{center}
$\sigma: c_{i}\to c'_{i}$, for $1\leq i\leq 2r$.
\end{center}

Then we have
\begin{center}
$B_{Q_{2}}(c_{i}^{\sigma}, c_{j}^{\sigma})=B_{Q_{1}}(c_{i}, c_{j})$ and $Q_{2}(c^{\sigma}_{i})=Q_{1}(c_{i})$,  for all $1\leq i, j\leq 2r$.
\end{center}

Thus $\sigma$ is an isometry of $Q_{1}$ and hence $Q_{1}$ and $Q_{2}$ have the same type. Therefore by Lemma~\ref{ch4lem32} either $Q_{2}$ is hyperbolic with $r\equiv 0$ or $1\Mod 4$, or $Q_{2}$ is elliptic with $r\equiv 2$ or $3\Mod 4$. \qed

Combining the results of this section, we obtain the following theorem.

\begin{guess}\label{ch4themsn}
Let $V$ be a vector space of dimension $2r$ over $\mathbb{F}_{2}$ with nondegenerate quadratic form $Q$. Then $V$ has a symmetric basis if and only if either $Q$ is hyperbolic and $r\equiv 0$ or $1 \Mod 4$, or $Q$ is elliptic and $r\equiv 2$ or $3 \Mod 4$.
\end{guess}

\section{A 2-Arc-Transitive Cover of Hypercubes}

Let $G$ be an extraspecial $2$-group of order $2^{2r+1}$ with identity {\bf 1} ($r\geqslant  1$). Let $Z=\langle z\rangle$ be the center of $G$. Then $Z\cong \mathbb{Z}_{2}$ and $G/Z\cong \mathbb{Z}_{2}^{2r}$ is elementary abelian. The commutator of any two elements in $G$ or the square of any element in $G$ lies in $Z$. So $G$ is a nilpotent group of class $2$. Define two functions $B: G/Z\times G/Z \to Z$ and $Q: G/Z\to Z$ as below: for any $Zx, Zy$ in $G/Z$, 
\begin{center}
$B(Zx, Zy)=[x,y]$,\\
$Q(Zx)=x^{2}$.
\end{center}
Then $Q$ is a quadratic form on $V=G/Z$ with associated bilinear form $B$. Note that if $B(Zx, Zy)={\bf 1}$ for some $x, y\in G$, then $x,y$ commute. So if $B(Zx, Zy)={\bf 1}$ for all $y\in G$, then  $Zx$ must be the identity in $G/Z$. Therefore $Q$ is nondegenerate on $G/Z$. Furthermore $B$ is symmetric  as $[y, x]=[x,y]^{-1}=[x, y]$ for all $x, y\in G$. 

We say that $G$ is an extraspecial $2$-group of {\it plus type} if $Q$ is hyperbolic, denoted by $2_{+}^{2r+1}$, and $G$ is an extraspecial $2$-group of {\it minus type} if $Q$ is elliptic, denoted by $2_{-}^{2r+1}$. It is known \cite{winter1972} that if $G=2_{+}^{2r+1}$, then it is the central product of $r$ dihedral groups $D_{8}$, otherwise $G$ is the central product of $r-1$ dihedral groups $D_{8}$ with one quaternion group $Q_{8}$. Also Winter proved that $\mathrm{Aut}(2_{\epsilon}^{2r+1})=2^{2r}. O^{\epsilon}(2r, 2)$ where $O^{\epsilon}(2r, 2)$ with $\epsilon\in\{+, -\}$ is an orthogonal group (see \cite[Theorem 1]{winter1972}).   

By Theorem~\ref{ch4themsn}, $G/Z$ has a symmetric basis if and only if either $G$ is of plus type with $r\equiv 0$ or $1 \Mod 4$, or $G$ is of minus type with $r\equiv 2$ or $3 \Mod 4$. Let $\mathcal{B}=\{Zg_{1}, Zg_{2}, \ldots, Zg_{2r}\}$ be a symmetric basis of $G/Z$, and so for all $1\leqslant  i\neq j\leqslant  2r$ we have that $Q(Zg_{i})=(g_{i})^{2}=\bf1$ and $B(Zg_{i}, Zg_{j})=[g_{i}, g_{j}]=z$. This implies that $g_{1}, g_{2}, \ldots, g_{2r}$ are involutions of $G$, and $g_{i}g_{j}=g_{j}g_{i}z$ for all distinct $i$ and $j$ with $1\leqslant  i, j\leqslant  2r$. Let $\Sigma=\mathrm{Cay}(G/Z, \mathcal{B})$ and $\Gamma=\mathrm{Cay}(G, S)$ with $S=\{g_{1}, g_{2}, \ldots, g_{2r}\}$. Thus $\Sigma$ is a $2r$-dimensional hypercube. 

To prove Theorem~\ref{ch4mainthem}, we first show that $\Gamma$ is a $2$-arc-transitive cover of $\Sigma$. The  {\it Frattini subgroup $\Phi(M)$} of a group $M$ is the intersection of all maximal subgroups of $M$. If $M$ is a $p$-group, then the {\it Frattini quotient $M/\Phi(M)$ of $M$} is isomorphic to $Z_{p}^{k}$ where $k$ is the smallest number of generators for $M$. Since $G$ is an extraspecial $2$-group, we have that $\Phi(G)=Z$. Since $G/Z\cong \mathbb{Z}_{2}^{2r}$ has a symmetric basis $\mathcal{B}$, by the Burnside Basis Theorem (see \cite[Theorem 11.12]{rose1994course}) we have that $G=\langle g_{1}, g_{2}, \ldots, g_{2r}\rangle$.

Let $g\in G$ where {\bf 1}$\neq g\neq z$. Then $Zg$ can be uniquely written as $Zg_{s_{1}}Zg_{s_{2}}\cdots Zg_{s_{t}}$ where $1\leqslant  s_{1}<s_{2}<\cdots <s_{t}\leqslant  2r$. So $g=z^{j}g_{s_{1}}\cdots g_{s_{t}}$ for $j=0$ or $1$.

For each $\sigma\in S_{2r}$, define a map $\tilde{\sigma}: G \to G$ by $z^{\tilde{\sigma}}=z$, and for all $g\in G$ with $g\neq z$, 
\begin{center}
$g^{\tilde{\sigma}}=z^{j}g_{s_{1}^{\sigma}}\cdots g_{s_{t}^{\sigma}}$
\end{center}
where $g=z^{j}g_{s_{1}}\cdots g_{s_{t}}$ for $j=0$ or $1$. Note that $(Zg_{1})(Zg_{2})=Zg_{1}g_{2}=Zg_{2}g_{1}$ for any $g_{1}$ and $g_{2}$ in $G$. In particular, for all $s_{i}, s_{k}$ we have $g_{s_{i}}g_{s_{k}}=z^{j}g_{s_{k}}g_{s_{i}}$  for some $j$ and so we can deduce that $\tilde{\sigma}$ is a homomorphism. 
Suppose that $h\in \mathrm{Ker}(\tilde{\sigma})$ and $h\neq{\bf 1}$. Since $h\neq z$, the element $h$ can be uniquely written as $z^{j}g_{s_{1}}\cdots g_{s_{t}}$ for some $s_{i}\in \{1, \ldots, 2r\}$. Since $h^{\tilde{\sigma}}=${\bf 1},  we have 
$$
Z=Z(h^{\tilde{\sigma}})=Zg_{s_{1}^{\sigma}}\cdots Zg_{s_{t}^{\sigma}}.
$$
This is a contradiction as $\{Zg_{1}, \ldots, Zg_{2r}\}$ is a basis of $G/Z$. Therefore, $\mathrm{Ker}(\tilde{\sigma})=${\bf 1}. Thus, $\tilde{\sigma}$ is injective and as $G$ is finite, it follows that $\tilde{\sigma}$ is surjective. Hence, $\tilde{\sigma}\in Aut(G)$. 

\begin{guess}\label{ch4them2}
$S_{2r}\leqslant  Aut(G)$.
\end{guess}

\pf Let $\phi: S_{2r}\to Aut(G)$ be the map defined as below: 
\begin{center}
$\phi: \sigma\to \tilde{\sigma}$, for each $\sigma\in S_{2r}$.
\end{center}
It is not hard to prove that $\phi$ is a homomorphism from $S_{2r}$ into $Aut(G)$. Let $K$ be the kernel of $\phi$, that is, 
$$
K=\{\sigma\in S_{2r}\mid \tilde{\sigma}=\bf1\}.
$$ 
Since $\tilde{\sigma}={\bf 1}$, we have that $g_{i}^{\tilde{\sigma}}=g_{i^{\sigma}}=g_{i}$ for all $1\leqslant  i\leqslant  2r$. Thus $\sigma=\bf1$, and so $K=\{{\bf 1}\}$. Therefore, $S_{2r}\leqslant  Aut(G)$. \qed
   
\begin{guess}\label{ch4them2arc}
The graph $\Gamma$ is a 2-arc-transitive Cayley graph of $G$.
\end{guess}

\pf Let $A=Aut(\Gamma)$. Then $S_{2r}\leqslant  A_{{\bf1}}$ as $S_{2r}$ fixes the identity element ${\bf1}$ of $G$. Let $N({\bf1})$ be the set of neighbours of ${\bf1}$ in $\Gamma$. Thus $S_{2r}$ is 2-transitive on $N({\bf 1})$. Since $\Gamma$ is vertex-transitive, we have that $\Gamma$ is $2$-arc-transitive.  \qed

\begin{guess}\label{ch4themcover}
$\Gamma$ is a 2-arc-transitive cover for an even-dimensional hypercube.
\end{guess}
\pf It follows from Theorem~\ref{ch4them2arc} that $G\rtimes S_{2r}\leqslant  \mathrm{Aut}(\Gamma)$ is 2-arc-transitive on $\Gamma$. Since $|G:Z(G)|\geqslant  4$, by \cite[Theorem 4.1]{praeger1993nan}  we have that  $\Gamma$ is a cover of $\Gamma_{Z(G)}\cong \Sigma$. \qed
\medskip

\medskip

Now we prove that $\Gamma$ is a normal Cayley graph for $G$. By a computation in {\textsc{Magma}} \cite{bosma1997magma}, we found that $\Gamma$ is a normal Cayley graph for $G$ for all $r\in\{1, 2, 3\}$. Next, we show that this is true for the general case, that is, for all $r\geqslant  4$.

Let $C$ be a cycle in $\Gamma$ with $V(C)=\{c_{1}, \ldots, c_{t_{c}}\}$ and $E=\{\{c_{t_{c}}, c_{1}\}\}\cup \{\{c_{i}, c_{i+1}\}\mid 1\leqslant  i\leqslant  t_{c}-1\}$ where $t_{c}$ is the length of $C$. There is a sequence $(s_{1}, \ldots, s_{t_{c}})$ induced by $C$ where $s_{i}\in S$ for $1\leqslant  i\leqslant  t_{c}$, such that $c_{1}=s_{t_{c}}c_{t_{c}}$ and $c_{i+1}=s_{i}c_{i}$ for $1\leqslant  i\leqslant  t_{c}-1$. Since  $s_{i}$ is an involution for all $1\leqslant i\leqslant 2r$, we have that 
\begin{equation}\label{ch4equ1}
s_{1}s_{2}\cdots s_{t_{c}}={\bf 1},
\end{equation}
and 
\begin{equation}\label{ch4equ4}
c_{i}= \begin{cases}
               s_{t_{c}}c_{t_{c}},               & i=1, \\
               (s_{i-1}s_{i-2}\cdots s_{1})c_{1},              &  2\leqslant  i\leqslant  t_{c}.
               \end{cases}
\end{equation}
The sequence is uniquely determined by $C$. Since $\Gamma$ is simple, we have the following lemma.

\begin{lemma}\label{ch4lem1} 
$s_{1}\neq s_{t_{c}}$, and $s_{i}\neq s_{i+1}$ for all $1\leqslant  i\leqslant  t_{c}-1$.
\end{lemma}
We call the sequence induced by $C$ the {\it cycle-sequence for $C$}. 

\medskip

\begin{lemma}\label{ch4lem2}
Suppose that $(s_{1}, \ldots, s_{n})$ is a sequence with $s_{i}\in S$ for all $1\leqslant  i\leqslant  n$ such that $s_{1}\cdots s_{n}\in Z$. Then for each $1\leqslant  i\leqslant  n$, the element $s_{i}$ appears an even number of times in the sequence. 
\end{lemma}

\pf Note that 
\begin{equation*}
Z=Zs_{1}\cdots s_{n}= (Zs_{1})\cdots(Zs_{n})=(Zs_{m_{1}})^{k_{1}}\cdots(Zs_{m_{t}})^{k_{t}}, 
\end{equation*}
where $\sum\limits_{1\leqslant  i\leqslant  t}k_{i}=n$ and $k_{i}$ is the number of $s_{m_{i}}$ in $(s_{1}, \ldots, s_{n})$. Since $\{Zs_{1}, \ldots, Zs_{2r}\}$ is a basis and $s_{i}$ is an involution for all $1\leqslant  i\leqslant  2r$, we have that $k_{i}$ is even for all $1\leqslant  i\leqslant  t$. \qed 

\medskip

\begin{figure}[h]
\begin{center}$
\begin{array}{c|c}
\includegraphics[scale=0.8]{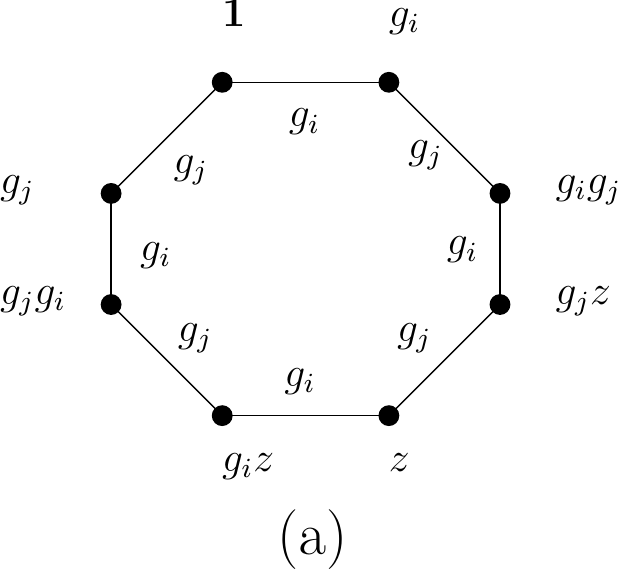} & \includegraphics[scale=0.8]{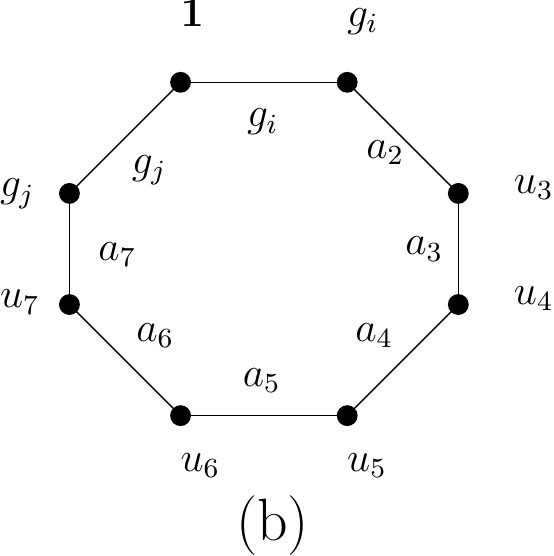} 
\end{array}$
\end{center}
\caption{The 8-cycle $C_{ij}$ and its image under $\rho\in A_{\bf1}^{[1]}$ ($1\leqslant  i, j\leqslant  2r$, $i\neq j$)}
\label{ch4fig1}
\end{figure}

Let $A_{\bf 1}^{[1]}$ be the automorphisms in $A_{\bf 1}$ that fix each vertex in $N({\bf 1})$. Recall that $N({\bf1})=S$. For each distinct $i$ and $j$ with $1\leqslant  i, j\leqslant  2r$, note that the sequence $(s_{1}, \ldots, s_{8})$ defined by 
$$
\hspace{180pt} s_{k}= \begin{cases}
               g_{i},               & \textit{if k is odd}, \\
               g_{j},              &  \textit{if k is even}.
               \end{cases}  \hspace{180pt} 
$$
is an 8-cycle in $\Gamma$ which we will denote by $C_{ij}$ (see Figure~\ref{ch4fig1}(a)). Let $c_{1}=\bf1$ and for each $2\leqslant  k\leqslant 8$, let $c_{k}=s_{k-1}s_{k-2}\cdots s_{1}$, in particular, we have that $c_{5}=g_{j}g_{i}g_{j}g_{i}=z$ and $c_{8}=g_{j}$.  

Let $\rho\in A_{\bf 1}^{[1]}$ and let $C_{ij}^{\rho}$ be the image of $C_{ij}$ under $\rho$ with $V(C_{ij}^{\rho})=\{u_{k}\mid u_{k}=c_{k}^{\rho}, 1\leqslant  k\leqslant  8\}$. Since $g_{i}, g_{j}\in N(1)$ and $\rho\in A_{\bf 1}^{[1]}$, we have that $u_{1}={\bf 1}$, $u_{2}=g_{i}$ and $u_{8}=g_{j}$. Let $\{a_{1}, \ldots, a_{8}\}$ be the cycle-sequence for $C_{ij}^{\rho}$, and so $a_{1}=g_{i}$ and $a_{8}=g_{j}$ (see Figure~\ref{ch4fig1}(b)).  

\begin{lemma}\label{ch4lem3}
Suppose that $z^{\rho}\neq z$. Then for $1\leqslant  i\leqslant  8$, the element $a_{i}$ appears exactly twice in the cycle sequence $(s_{1}, \ldots, s_{8})$. Moreover, $a_{1}, a_{2}, a_{3}, a_{4}$ are pairwise distinct, and $a_{5}, a_{6}, a_{7}, a_{8}$ are pairwise distinct.
\end{lemma}

\pf By Lemma~\ref{ch4lem1} we have that $a_{1}\neq a_{8}$ and $a_{i}\neq a_{i+1}$ for all $1\leqslant  i\leqslant  7$. Suppose that $a_{3}=a_{1}$. Then by (\ref{ch4equ4}), we have that 
\begin{align*}
u_{5} &= a_{4}a_{3}a_{2}a_{1}\\
& = a_{4}a_{1}a_{2}a_{1}\\
& = a_{4}a_{2}z\\
& = a_{2}a_{4}.
\end{align*}
If $a_{2}=a_{4}$, then $u_{5}={\bf 1}$, which is a contradiction. Thus $a_{2}\neq a_{4}$, which implies that there is a $2$-path in $\Gamma$ connecting {\bf 1} and $u_{5}$. Since $\rho\in A$ and $z^{\rho}=u_{5}$, we have that there is a $2$-path in $\Gamma$ connecting {\bf 1} and $z$, which leads to a contradiction. Hence $a_{3}\neq a_{1}$. If $a_{4}=a_{1}$, then 
\begin{align*}
u_{5} &= a_{4}a_{3}a_{2}a_{1}\\
& = a_{1}a_{3}a_{2}a_{1}\\
& = a_{3}a_{2},
\end{align*}
and by the same arguments, we deduce that $a_{4}\neq a_{1}$. Suppose that $a_{2}=a_{4}$. Then $u_{5} = a_{4}a_{3}a_{2}a_{1}=a_{1}a_{3}$, which leads to a contradiction by the same arguments. Thus $a_{1}, a_{2}, a_{3}, a_{4}$ are pairwise distinct. By (\ref{ch4equ4}) we have that 
\begin{equation*}
u_{8}=a_{7}a_{6}a_{5}u_{5}.
\end{equation*}
Since ${\bf 1}=u_{1}=a_{8}u_{8}$, we have that ${\bf 1}=a_{8}a_{7}a_{6}a_{5}u_{5}$. Since $S$ consists of involutions, we have that $u_{5}=a_{5}a_{6}a_{7}a_{8}$. Then by the same arguments, we may conclude that $a_{5}, a_{6}, a_{7}, a_{8}$ are pairwise distinct. Therefore, by Lemma~\ref{ch4lem2}, we conclude that each term in $\{a_{1}, \ldots, a_{8}\}$ appears exactly twice. \qed

\medskip

\begin{lemma}\label{ch4lem4}
The cycle $C_{ij}$ is fixed pointwise by $\rho$.
\end{lemma}

\pf  We first show that $z^{\rho}=z$. Suppose to the contrary that $z^{\rho}\neq z$. Since $u_{5}=z^{\rho}$, we have that $(N(z)\cap C_{ij})^{\rho}=\{g_{i}z, g_{j}z\}^{\rho}=\{u_{4}, u_{6}\}$. Recall that $u_{4}=a_{3}a_{2}a_{1}$ and $u_{6}=a_{5}a_{4}a_{3}a_{2}a_{1}$. 

Suppose that $(g_{i}z)^{\rho}=u_{4}$. Let $g_{k}\in S$ be such that $g_{k}\neq a_{i}$ for all $1\leqslant  i\leqslant  3$, and let $C_{ik}$ be an $8$-cycle. Note that $C_{ik}$ has the same shape as in Figure~\ref{ch4fig1}(a). Thus $g_{i}z$ and $g_{k}$ are connected by a $2$-path in $C_{ik}$, and so we have that $(g_{i}z)^{\rho}$ and $g_{k}^{\rho}$ are connected by a $2$-path, that is, $u_{4}$ and $g_{k}$ are connected by a $2$-path in $\Gamma$. We may assume that $x_{2}x_{1}u_{4}=g_{k}$ for some $x_{1}, x_{2}\in S$ where $x_{1}\neq x_{2}$, and so $x_{2}x_{1}a_{3}a_{2}a_{1}g_{k}={\bf 1}$. By Lemma~\ref{ch4lem2} and Lemma~\ref{ch4lem3}, we must have that $g_{k}=a_{i}$ for some $1\leqslant  i\leqslant  3$, which leads to a contradiction. Thus $(g_{i}z)^{\rho}\neq u_{4}$, and so $(g_{i}z)^{\rho}=u_{6}$. 

By Lemma~\ref{ch4lem3}, $a_{5}=a_{j}$ for some $1\leqslant  j\leqslant  3$. Suppose that $a_{5}=a_{1}$. Then $(g_{i}z)^{\rho}=u_{6}=a_{1}a_{4}a_{3}a_{2}a_{1}=a_{4}a_{3}a_{2}z$. Let $g_{k'}\in S$ such that $g_{k'}\neq a_{i}$ for $2\leqslant  i\leqslant  4$, and let $C_{ik'}$ be the corresponding $8$-cycle. Since $g_{i}z$ and $g_{k'}$ are joined by a 2-path, we have that $u_{6}$ and $g_{k'}$ are joined by a 2-path, and so we have $y_{2}y_{1}a_{4}a_{3}a_{2}g_{k'}=z$ for some $y_{1}, y_{2}\in S$ where $y_{1}\neq y_{2}$. Then by similar arguments,  we conclude that $g_{k'}=a_{j}$ for some $2\leqslant  j\leqslant  4$, which is a contradiction. For the remaining two cases where $a_{5}=a_{2}$ or $a_{5}=a_{3}$, we can obtain contradictions following similar arguments. Hence $(g_{i}z)^{\rho}\neq u_{6}$, which leads to a contradiction to the fact that $(g_{i}z)^{\rho}\in \{u_{4}, u_{6}\}$. Therefore, $z^{\rho}=z$, that is, $u_{5}=z$. 

Recall that $u_{5}=a_{4}a_{3}a_{2}a_{1}=a_{5}a_{6}a_{7}a_{8}$ where $a_{1}=g_{i}$ and $a_{8}=g_{j}$. By Lemma~\ref{ch4lem2} we have that $a_{1}=a_{3}=g_{i}$, $a_{2}=a_{4}$, $a_{5}=a_{7}$ and $a_{6}=a_{8}=g_{j}$ as $u_{5}=z$ and $a_{i}\neq a_{i+1}$ for all $1\leqslant  i\leqslant  7$. If $a_{2}=a_{4}=g_{j}$ and $a_{5}=a_{7}=g_{i}$, then $C_{ij}$ and $C_{ij}^{\rho}$ have the same cycle sequence, that is, $C_{ij}^{\rho}=C_{ij}$. 

Suppose to the contrary that $a_{2}=a_{4}\neq g_{j}$. Since $u_{5}=z$, we have that $u_{4}=a_{4}z=a_{2}z$ and $u_{6}=a_{5}z$. Recall that $\{g_{i}z, g_{j}z\}^{\rho}=\{u_{4}, u_{6}\}$. Since $a_{2}\in S$, we may assume that there exists $1\leqslant  t\leqslant  2r$ such that $a_{2}=g_{t}$ with $t\neq i$ and $t\neq j$. Suppose that $(g_{i}z)^{\rho}=u_{4}=a_{2}z$, that is, $(g_{i}z)^{\rho}=g_{t}z$. Let $C_{it}$ be the corresponding 8-cycle. Thus $g_{t}$ and $g_{i}z$ are connected by a 2-path in $C_{it}$, and so we have that $g_{t}$ and $(g_{i}z)^{\rho}$ are connected by a 2-path. Thus there exists $w_{1}, w_{2}\in S$ such that $w_{1}\neq w_{2}$ and $w_{2}w_{1}g_{t}=g_{t}z$, that is, $w_{1}w_{2}=z$, which is a contradiction as $Zw_{1}$ and $Zw_{2}$ are base elements. Hence $(g_{i}z)^{\rho}=u_{6}=a_{5}z$, that is, $(g_{j}z)^{\rho}=u_{4}=a_{2}z$.

Recall that $a_{2}=g_{t}$ where $t\neq i$ and $t\neq j$. Let $C_{jt}$ be the corresponding 8-cycle. Note that $C_{jt}$ has the same shape as in Figure~\ref{ch4fig1}(a). Thus $g_{j}z$ and $g_{t}$ are connected by a 2-path. Hence $(g_{j}z)^{\rho}=a_{2}z=g_{t}z$ and $g_{t}$ are connected by a 2-path. Suppose that there exist $u_{1}, u_{2}\in S$ such that $u_{2}u_{1}g_{t}=g_{t}z$. Thus we have that $u_{1}u_{2}=z$, which is a contradiction. Hence $a_{2}=a_{4}=g_{j}$. 

By a similar argument, we may conclude that $a_{5}=a_{7}=g_{i}$. Therefore $C_{ij}^{\rho}=C_{ij}$. \qed

\medskip

\begin{corollary}\label{ch4coro1}
Let $v\in N({\bf 1})$. Then $\rho\in A_{v}^{[1]}$. 
\end{corollary}

\medskip

\pf  Let $w\in N(v)$. We may assume that $v=g_{1}$, and so $w=g_{1}g_{i}$ for some $1\leqslant  i\leqslant  2r$. If $i=1$, then $w={\bf 1}$, and so $\rho$ fixes ${\bf 1}$. Suppose that $i\neq 1$, and consider the $8$-cycle $C_{1i}$. By Lemma~\ref{ch4lem4}, $C_{1i}^{\rho}=C_{1i}$, that is, $\rho$ fixes each vertex on $C_{1i}$, and so $w^{\rho}=w$. Therefore $\rho\in A_{v}^{[1]}$. \qed

\begin{lemma}\label{ch4lem5}
$A_{\bf 1}^{[1]}=\{\bf 1\}$. 
\end{lemma}

\pf By Corollary~\ref{ch4coro1} we have that $A_{\bf1}^{[1]}\leqslant  A_{v}^{[1]}$. Since $\Gamma$ is vertex-transitive and finite, we have that $|A_{\bf1}^{[1]}|=|A_{v}^{[1]}|$, which implies that $A_{\bf1}^{[1]}=A_{v}^{[1]}$. Since $\Gamma$ is vertex-transitive, it follows that $A_{u}^{[1]}=A_{w}^{[1]}$ for each $u\in V(\Gamma)$ and $w\in N(u)$. Thus by connectivity $A_{\bf1}^{[1]}=A_{u}^{[1]}$ for all $u\in V(\Gamma)$, and so $A_{\bf1}^{[1]}=\{\bf 1\}$. \qed

\medskip

{\it\small Proof of Theorem~\ref{ch4mainthem}.}  By Lemma~\ref{ch4lem5} we have $A_{\bf 1}\cong A_{\bf 1}^{N(\bf 1)}=S_{2r}$. Further by Theorem~\ref{ch4them2} we obtain $A_{\bf 1}=Aut(G, S)$. Therefore $A=G\rtimes S_{2r}$, namely, $\Gamma$ is a normal Cayley graph for $G$. Also in Theorem~\ref{ch4themcover} we have proved that $\Gamma$ is a $2$-arc-transitive cover of a hypercube of dimension $2r$. This completes the proof of Theorem~\ref{ch4mainthem}. \qed

\section{Acknowledgements}

This research is part of my Ph.D project at the Centre for the Mathematics of Symmetry and Computation
of the University of Western Australia. I would like to express my gratitude to my supervisors Dr. Michael Giudici and Dr. Cai Heng Li, who offered continuous advice and encouragement throughout the progress of this paper. 

I would also like to thank the referees of this paper for their great suggestions, which helped me a lot in improving the presentation of this paper.

\bibliographystyle{acm}
\bibliography{ref}

\end{document}